\documentstyle[12pt]{amsart}

\newtheorem{th}{Theorem}[section] 
\newtheorem{lem}[th]{Lemma} 
\newtheorem{cor}[th]{Corollary} \newtheorem{prp}[th]{Proposition}

\theoremstyle{definition} 
\newtheorem{dfn}[th]{Definition} 
 \newtheorem{conj}[th]{Conjecture}
 \newtheorem{question}[th]{Question}

\theoremstyle{remark} 

\newcommand{\rar}{\rightarrow}

\newcommand{\bfz}{{\Bbb{Z}}} 
\newcommand{\bfr}{{\Bbb{R}}}

\newcommand{\QED}{\qed\\}
\newcommand{\LA}{\langle}
\newcommand{\RA}{\rangle}

\newcommand{\Span}{{\operatorname{Span}}}

\begin{document}

\noindent%

\title[Semistable reduction]{Semistable reduction in characteristic 0
for families of surfaces and three-folds}
\author{K. Karu}
\address{Department of Mathematics\\ Boston University\\ 111
Cummington\\ Boston, MA 02215, USA}
\email{kllkr@@math.bu.edu}
\date{\today}

\maketitle

\section{Introduction}
In \cite{ak} the semistable reduction of a morphism $F:X\rar B$ was
stated as a problem in the combinatorics of polyhedral complexes. In
this
paper we solve it in the case 
when the relative dimension of $F$ is no bigger than three.  

First we   recall the setup of the problem from \cite{ak}. The ground
field
$k$ will be algebraically closed of characteristic zero. 

\begin{dfn} A flat morphism $F:X\rar B$ of nonsingular projective
varieties
is semistable if in local analytic coordinates $x_1,\ldots,x_n$ at $x\in
X$
and $t_1,\ldots,t_m$ at $b\in B$ the morphism $F$ is given by
\[ t_i = \prod_{j=l_{i-1}+1}^{l_i} x_j \]
where $0=l_0<l_1<\ldots<l_m\leq n$.
\end{dfn}

The conjecture of semistable reduction states that

\begin{conj}\label{conj-ssr} Let $F: X\rar B$ be a surjective morphism
with
geometrically integral generic fiber. There exist an alteration (proper
surjective generically finite 
morphism) $B'\rar B$ and a modification (proper biratonal morphism)
$X'\rar
X\times_B B'$ such that $X'\rar B'$ is semistable.
\end{conj}

Conjecture~\ref{conj-ssr} was proved in \cite{te} (main theorem of
Chapter~2) in 
case when $B$ is a curve. A weak version of the conjecture was proved in
\cite{ak} for arbitrary $X$ and $B$. In both cases the proof proceeds by
reducing $F$ to a morphism of toroidal embeddings, stating the problem
in
terms of the associated polyhedral complexes, and solving the
combinatorial
problem. 

\subsection{Polyhedral complexes}
We consider (rational, conical) polyhedral complexes
$\Delta=(|\Delta|,\{\sigma\},\{N_\sigma\})$ consisting of a collection
of
lattices $N_\sigma\cong\bfz^n$ and rational full cones $\sigma\subset
N_\sigma\otimes\bfr$ with a vertex. The cones $\sigma$ are glued
together to
form the space $|\Delta|$ so that the usual axioms of polyhedral
complexes
hold: 
\begin{enumerate}
\item If $\sigma\in\Delta$ is a cone, then every face $\sigma'$ of
$\sigma$
is also in $\Delta$, and $N_{\sigma'}=N_\sigma|_{\Span(\sigma')}$.
\item The intersection of two cones $\sigma_1\cap\sigma_2$ is a
face of both of them, $N_{\sigma_1\cap\sigma_2} =
N_{\sigma_1}|_{\Span(\sigma_1\cap\sigma_2)} =
N_{\sigma_2}|_{\Span(\sigma_1\cap\sigma_2)}$. 
\end{enumerate}
A morphism $f:\Delta_X\rar\Delta_B$ of polyhedral complexes
$\Delta_X=(|\Delta_X|,\{\sigma\},\{N_\sigma\})$ and
$\Delta_B=(|\Delta_B|,\{\tau\},\{N_\tau\})$ is a compatible collection
of
linear maps $f_\sigma: (\sigma,N_\sigma)\rar(\tau,N_\tau)$; i.e. if
$\sigma'$
is a face of $\sigma$ then $f_{\sigma'}$ is the restriction of
$f_\sigma$. We
will only consider morphisms $f:\Delta_X\rar\Delta_B$ such that
$f_\sigma^{-1}(0)\cap\sigma=\{0\}$ for all $\sigma\in\Delta_X$.

\begin{dfn} A surjective morphism $f:\Delta_X\rar\Delta_B$ such that
$f^{-1}(0)=\{0\}$ is semistable if
\begin{enumerate}
\item $\Delta_X$ and $\Delta_B$ are nonsingular.
\item For any cone $\sigma\in\Delta_X$, we have $f(\sigma)\in\Delta_B$
and
$f(N_\sigma)=N_{f(\sigma)}.$ 
\end{enumerate}
We say that $f$ is weakly semistable if it satisfies the two properties
except that $\Delta_X$ may be singular.
\end{dfn}

The following two operations are allowed on $\Delta_X$ and $\Delta_B$:
\begin{enumerate}
\item Projective subdivisions $\Delta_X'$ of $\Delta_X$ and $\Delta_B'$
of
$\Delta_B$ such that $f$ induces a morphism 
$f':\Delta_X'\rar\Delta_B'$;
\item Lattice alterations: let
$\Delta_X'=(|\Delta_X|,\{\sigma\},\{N_\sigma'\}),
\Delta_B'=(|\Delta_B|,\{\tau\},\{N_\tau'\})$, for some compatible
collection
of sublattices $N_\tau'\subset N_\tau$, $N_\sigma'=f^{-1}(N'_\tau)\cap
N_\sigma$, and let $f':\Delta_X'\rar\Delta_B'$ be the morphism induced
by
$f$.
\end{enumerate}

\begin{conj}\label{main-conj} Given a surjective morphism
$f:\Delta_X\rar\Delta_B$, such that $f^{-1}(0)=\{0\}$, there exists a
projective subdivision $f':\Delta_X'\rar\Delta_B'$ followed by a lattice
alteration $f'':\Delta_X''\rar\Delta_B''$  so that $f''$ is semistable. 
\[
\begin{array}{lclcl}
\Delta_{X''} & \rar & \Delta_{X'} & \rar & \Delta_{X} \\
\downarrow f'' & & \downarrow f' & & \downarrow f \\
\Delta_{B''} & \rar & \Delta_{B'} & \rar & \Delta_{B}
\end{array}
\]
\end{conj}

The importance of Conjecture~\ref{main-conj} lies in the fact that
it implies Conjecture~\ref{conj-ssr}
(Proposition~8.5 in \cite{ak}). In the case when $\dim(\Delta_B)=1$,
Conjecture~\ref{main-conj} was proved in \cite{te} (main theorem of
Chapter~3). In \cite{ak} (Theorem~0.3) the conjecture was proved with
semistable
replaced by weakly semistable. The main result of this paper is

\begin{th}\label{main-thm} Conjecture~\ref{main-conj} is true if $f$ has
relative dimension $\leq 3$. Hence, Conjecture~\ref{conj-ssr} is true if
$F$
has relative dimension $\leq 3$.
\end{th}

The relative dimension of a linear map $f:\sigma\rar\tau$ of cones
$\sigma,
\tau$ is $\dim(\sigma)-\dim(f(\sigma))$. The relative dimension of
$f:\Delta_X\rar\Delta_B$ is by definition the maximum of the relative
dimensions of $f_\sigma:\sigma\rar\tau$ over all $\sigma\in\Delta_X$. 
If $F:X\rar B$ is a morphism of toroidal embeddings of relative
dimension $d$,
then the associated morphism of polyhedral complexes
$f:\Delta_X\rar\Delta_B$ has relative dimension $\leq d$ because in
local models the relative dimension of $F$ is no bigger than the rank of
the
kernel of $f: N_\sigma\rar N_\tau$. Thus, the second statement of the
theorem
follows from the  first.

\subsection{Notation} We will use notations from \cite{te} and
\cite{fu}. For a cone $\sigma\in N\otimes\bfr$ we write
$\sigma={\LA}v_1,\ldots,v_n{\RA}$ if $v_1,\ldots,v_n$ lie on the
1-dimensional edges
of $\sigma$ and generate it. If $v_i$ are the first lattice points along
the
edges we call them primitive points of $\sigma$. For a simplicial cone
$\sigma$ with primitive points $v_1,\ldots,v_n$, the multiplicity of
$\sigma$
is  
\[ m(\sigma,N_\sigma) = [N_\sigma:\bfz v_1\oplus\ldots\oplus\bfz
v_n]. \] 
A polyhedral complex $\Delta$ is nonsingular if and only if
$m(\sigma,N_\sigma)=1$ for all $\sigma\in\Delta$. To compute the
multiplicity
of $\sigma$ we can count the representatives $w\in N_\sigma$ of classes
of
$N_\sigma/\bfz v_1\oplus\ldots\oplus\bfz v_n$ in the form
\[ w=\sum_{i}\alpha_i v_i, \qquad 0\leq\alpha_i<1.\]
Such points $w$ were called Waterman points of $\sigma$ in \cite{te}.
Also notice that since the multiplicity of a
face of $\sigma$ is no bigger than the multiplicity of $\sigma$, to
compute
the multiplicity of $\Delta$ it suffices to consider maximal cones only.

If $\Delta_X$ and $\Delta_B$ are simplicial, we we say that
$f:\Delta_X\rar\Delta_B$ is simplicial if $f(\sigma)\in\Delta_B$ for all
$\sigma\in\Delta_X$. Assume that $f$ is simplicial. Let
$u_1,\ldots,u_n$ be the primitive points of $\Delta_B$, and
$m_1,\ldots,m_n$
positive integers. By taking the $(m_1,\ldots,m_n)$ sublattice at
$u_1,\dots,u_n$ we mean the lattice alteration $N_\tau'=\bfz[m_{i_1}
u_{i_1},\ldots, m_{i_l}u_{i_l}]$ where $\tau\in\Delta_B$ has primitive
points
$u_{i_1},\dots,u_{i_l}$.

For cones $\sigma_1,\sigma_2\in\Delta$ we write $\sigma_1\leq\sigma_2$
if
$\sigma_1$ is a face of $\sigma_2$. 

A subdivision $\Delta'$ of $\Delta$ is called projective if there exists
a
homogeneous piecewise linear function $\psi:|\Delta|\rar\bfr$ taking
rational
values on 
the lattice points (a good function for short) such that the maximal
cones of
$\Delta'$ are exactly the maximal pieces in which $\psi$ is linear.   

\subsection{Acknowledgment} The suggestion to write up the proof of
semistable reduction for low relative dimensions came from Dan
Abramovich.

\section{Joins}
For cones $\sigma_1,\sigma_2\in\bfr^N$ lying in complementary planes: 
$\mbox{Span}(\sigma_1)\cap\mbox{Span}(\sigma_2)=\{0\}$, the join of
$\sigma_1$
and $\sigma_2$ is 
$\sigma_1*\sigma_2=\sigma_1+\sigma_2$. 
Let $\sigma$ be a simplicial cone $\sigma=\sigma_1*\ldots*\sigma_n$. If
$\sigma_i'$ is a subdivision of $\sigma_i$ for all $i=1,\ldots,n$, we
define
the join
\[ \sigma' = \sigma_1' * \ldots * \sigma_n' \]
as the set of cones $\rho = \rho_1+\ldots+\rho_n$, where
$\rho_i\in\sigma_i'$. 

Let $f:\Delta_X\rar\Delta_B$ be a simplicial map of simplicial
complexes. For $u_i$ a primitive point of $\Delta_B$, $i=1,\ldots,n$,
let
$\Delta_{X,i}=f^{-1}(\bfr_+ u_i)$ be the simplicial subcomplex of
$\Delta_X$. If $\Delta_{X,i}'$ is a subdivision of $\Delta_{X,i}$ for
$i=1,\ldots,n$, we can define the join
\[ \Delta_X' = \Delta'_{X,1} *\ldots*\Delta'_{X,n} \]
by taking joins inside all cones $\sigma\in\Delta_X$. This is well
defined by
the assumption that $f^{-1}(0)=\{0\}$.

\begin{lem} If $\Delta_{X,i}'$ are projective subdivisions of
$\Delta_{X,i}$
then the join $\Delta_{X}'$ is a projective subdivision of $\Delta_{X}$.
\end{lem}

{\bf Proof.} Let $\psi_i$ be good functions for $|\Delta_{X,i}'|$.
Extend
$\psi_i$ linearly to the entire $|\Delta_X'|$ by setting
$\psi_i(|\Delta_{X,j}'|)=0$ for $j\neq i$. Clearly, $\psi = \sum_i
\psi_i$
is a good function defining the subdivision $\Delta_{X}'$.
\QED

Consider $f|_{\Delta_{X,i}}: \Delta_{X,i}\rar\bfr_+ u_i$. By the main
theorem
of Chapter~2 in \cite{te} there 
exist a subdivision  $\Delta_{X,i}'$ of $\Delta_{X,i}$ and an
$m_i\in\bfz$
such that after taking the $m_i$-sublattice at $u_i$ we have
$f'|_{\Delta_{X,i}'}$ semistable. Now let 
$\Delta_X'$ be the join of $\Delta_{X,i}'$, and take the
$(m_1,\ldots,m_n)$-sublattice at $(u_1,\ldots,u_n)$. Then
$f':\Delta_X'\rar\Delta_B'$ is a simplicial map and
$f'|_{\Delta_{X,i}'}$ is 
semistable. 

We can also see that the multiplicity of $\Delta_X'$ is not bigger than
the
multiplicity of $\Delta_X$. Let $\sigma\in\Delta_X$ have primitive
points
$v_i$ and let $\sigma'\subset\sigma$ be a maximal cone in the
subdivision
with primitive points $v_i'$. The multiplicity of $\sigma'$ is the
number of Waterman points $w'\in N_\sigma'$
\[ w'=\sum_{i} \alpha_{i} v_i', \qquad 0\leq\alpha_{i}<1.\]
We show that the set of Waterman points of $\sigma'$ can be mapped
injectively into the set of Waterman points of $\sigma$, hence the
multiplicity of $\sigma'$ is not bigger than the multiplicity of
$\sigma$. Write 
\[ w'=\sum_{i}(\beta_i+b_i) v_i, \qquad 0\leq\beta_{i}<1, \qquad
b_i\in\bfz_+.\] 
Then $w=\sum_{i}\beta_i v_i \in N_\sigma$ is a Waterman point of
$\sigma$. If
different $w_1', w_2'$ give the same $w$, then $w_1'-w_2' \in
N_\sigma'\cap\bfz\{v_i\} = \bfz\{v_i'\}$, hence $w_1'-w_2'=0$. 

\section{Modified barycentric subdivisions}
Let $f:\Delta_X\rar\Delta_B$ be a simplicial morphism of simplicial
complexes. Consider the 
barycentric subdivision $BS(\Delta_B)$ of $\Delta_B$. The
1-dimensional cones of $BS(\Delta_B)$ are $\bfr_+\hat{\tau}$  where
$\hat{\tau}=\sum u_i$ is the barycenter of a cone 
$\tau\in\Delta_B$ with primitive points $u_1,\ldots,u_m$. A cone
$\tau'\in
BS(\Delta_B)$ is spanned by $\hat{\tau}_1,\ldots,\hat{\tau}_k$, where
$\tau_1\leq\tau_2\leq\ldots\leq\tau_k$ is a chain of cones in
$\Delta_B$.

In general, $f$ does not induce a morphism $BS(\Delta_X)\rar
BS(\Delta_B)$. For that we need to modify the barycenters
$\hat{\sigma}$ of cones $\sigma\in\Delta_X$.

\begin{dfn} The data of {\bf modified barycenters} consists of
\begin{enumerate}
\item A subset of cones $\tilde{\Delta}_X\subset\Delta_X$.
\item For each $\sigma\in\tilde{\Delta}_X$ a point $b_\sigma\in
\mbox{int}(\sigma)\cap N_\sigma$ such that
$f(b_\sigma)\in\bfr_+\hat{\tau}$
for some $\tau\in\Delta_B$.
\end{enumerate}
\end{dfn}

Recall that for any total order $\prec$ on the set of cones in
$\Delta_X$
refining the partial order $\leq$, the
barycentric subdivision $BS(\Delta_X)$ can be realized as a
sequence of star subdivisions at the barycenters $\hat{\sigma}$ of
$\sigma\in\Delta_X$ in the descending order $\prec$.  

\begin{dfn} Given modified barycenters $(\tilde{\Delta}_X,\{b_\sigma\})$
and
a total order $\prec$ on $\Delta_X$ refining the partial order
$\leq$, the {\bf modified barycentric subdivision}
$MBS_{\tilde{\Delta}_X,\{b_\sigma\},\prec}(\Delta_X)$ is the sequence of
star
subdivisions at $b_\sigma$ for $\sigma\in\tilde{\Delta}_X$ in the
descending
order $\prec$. 
\end{dfn}

To simplify notations, we will write $MBS(\Delta_X)$ instead of
$MBS_{\tilde{\Delta}_X,\{b_\sigma\},\prec}(\Delta_X)$. 
By definition, $MBS(\Delta_X)$ is a projective simplicial subdivision of
$\Delta_X$. As in the case of the ordinary barycentric subdivision, the
cones of
$MBS(\Delta_X)$  can be characterized by chains of cones in $\Delta_X$.
We
may assume that the 1-dimensional cones of $\Delta_X$ are all in
$\tilde{\Delta}_X$. For a cone $\sigma\in\Delta_X$ let $\tilde{\sigma}$
be
the maximal face of $\sigma$ (w.r.t. $\prec$) in $\tilde{\Delta}_X$.
Given a
chain of cones $\sigma_1\leq\ldots\leq\sigma_k$ in 
$\Delta_X$, the cone spanned by $b_{\tilde{\sigma}_1}, \ldots,
b_{\tilde{\sigma}_k}$ is a subcone of $\sigma_k$. Let $C(\Delta_X)$ be
the set
of all such cones corresponding to chains
$\sigma_1\leq\ldots\leq\sigma_k$ in
$\Delta_X$. 

\begin{prp} $C(\Delta_X)=MBS(\Delta_X)$.
\end{prp}

{\bf Proof.} Let $BS(\Delta_X)$ be the ordinary barycentric subdivision
of
$\Delta_X$. Both $C(\Delta_X)$ and $MBS(\Delta_X)$ are obtained from
$BS(\Delta_X)$ by moving the barycenters $\hat{\sigma}$ (and everything
attached to them) to the new position $b_{\tilde{\sigma}}$ for all
$\sigma\in\Delta_X$ in the descending order $\prec$.  
\QED

\begin{cor}\label{cor-simpl} If $f(\tilde{\sigma})=f(\sigma)$ for all
$\sigma\in\Delta_X$ then 
$f$ induces a simplicial map $f':MBS(\Delta_X)\rar BS(\Delta_B)$.
\end{cor}

{\bf Proof.} 
Let $\sigma'\in MBS(\Delta_X)$ correspond to a chain
$\sigma_1\leq\ldots\leq\sigma_k$. Then we have a chain of cones
$f(\sigma_1)\leq
\ldots\leq f(\sigma_k)$ in $\Delta_B$. The assumption that
$f(\tilde{\sigma}_i)=f(\sigma_i)$ implies that
$f(b_{\tilde{\sigma}_i})\in\bfr_+
\widehat{f}(\sigma_i)$, hence the cone 
${\LA}b_{\tilde{\sigma}_1},\ldots,b_{\tilde{\sigma}_k}{\RA}$ maps onto
the cone 
${\LA}\widehat{f}(\sigma_1), \ldots, \widehat{f}(\sigma_k){\RA}\in
BS(\Delta_B)$.
\QED

The hypothesis of the corollary is satisfied, for example, if for any
$\sigma\in\Delta_X$  with $f(\sigma)=\tau\in\Delta_B$ and for any face
$\sigma_1\leq\sigma$ such that $\sigma_1\in\tilde{\Delta}_X$,
$f(\sigma_1)\neq\tau$, there exists $\sigma_2\in\tilde{\Delta}_X$ such
that
$\sigma_1\leq\sigma_2\leq\sigma$ and $f(\sigma_2)=\tau$:
\[
\begin{array}{ccccc}
\sigma_1 & \leq & \sigma_2 & \leq & \sigma \\
\downarrow & & \downarrow & & \downarrow \\
\tau_1 & \leq & \tau & = & \tau \\
\end{array}
\]
Indeed, $\tilde{\sigma}\neq\sigma_1$ because $\sigma_1\prec\sigma_2$.

\subsection{Example} \label{example}
Assume that $f:\Delta_X\rar\Delta_B$ is a simplicial map of simplicial
complexes taking primitive
points of $\Delta_X$ to primitive points of $\Delta_B$ (e.g. $\Delta_X$
is
simplicial and $f$ is weakly semistable). Then for a cone
$\sigma\in\Delta_X$
such that $f:\sigma\stackrel{\simeq}{\rar}\tau$, we have
$f(\hat{\sigma})=\hat{\tau}$. 

Let $\tilde{\Delta}_X=\bar{\Delta}_X = \{\sigma\in\Delta_X: f|_\sigma
\mbox{is
injective}\}$, $b_\sigma=\hat{\sigma}$. In this case $\tilde{\sigma}$ is
the
maximal face of $\sigma$ (w.r.t. $\prec$) such that $f|_\sigma$ is
injective. Clearly, the hypothesis of the lemma is satisfied, and we
have a
simplicial map $f':MBS(\Delta_X)\rar BS(\Delta_B)$.

Next we compute the multiplicity of $MBS(\Delta_X)$. Let
$\sigma\in\Delta_X$
have primitive points $v_1,\ldots,v_n$, and let $\sigma'\subset\sigma$
be a
maximal cone in the subdivision, corresponding to the chain 
\[ \LA v_1\RA \leq \LA v_1,v_2\RA\leq\ldots\leq \LA v_1,\ldots,v_n\RA.
\]
Since $\tilde{\rho}\subset\rho$ for any $\rho$, the primitive points of
$\sigma'$ can be written as
\[
\begin{array}{lllllll}
v_1' &=& a_{11} v_1 & & & & \\
v_2' &=& a_{21} v_1 & + & a_{22} v_2 & & \\
     & \cdots & & & & & \\
v_n' &=& a_{n 1} v_1 & + & \ldots & + & a_{n n} v_n
\end{array}
\]
for some $0\leq a_{i j}$. The multiplicity of $\sigma'$ is $a_{1 1}\cdot
a_{2
2} \cdots a_{n n}$ times the multiplicity of $\sigma$. In case when
$b_\rho$
are barycenters $\hat{\rho}$, all $a_{i j}\leq 1$, hence the
multiplicity of
$\sigma'$ is not bigger than the multiplicity of $\sigma$.

\section{Reducing the multiplicity of $\Delta_X$}

Let $f:\Delta_X\rar\Delta_B$ be weakly semistable and $\Delta_X$
simplicial
(i.e. $\Delta_B$ is nonsingular, $\Delta_X$ is simplicial, and $f$ is a
simplicial map taking primitive points of $\Delta_X$ to primitive points
of
$\Delta_B$). Notice that if $\bar{\Delta}_X$ is as in
Example~\ref{example},
then $\bar{\Delta}_X$ is nonsingular, and
$f(\hat{\sigma})=\widehat{f}(\sigma)$ for any $\sigma\in\bar{\Delta}_X$.

A singular simplicial cone $\sigma\in\Delta_X$ with primitive points
$v_1,\ldots,v_n$ contains a Waterman point $w\in N_\sigma$,   
\[ w=\sum_{i} \alpha_{i} v_i, \qquad 0\leq\alpha_{i}<1, \qquad
\sum_i\alpha_i>0.\]
The star subdivision of $\sigma$ at $w$
has multiplicity strictly less than the multiplicity of $\sigma$.

We will show in this section that if every singular cone of $\Delta_X$
contains a Waterman point $w$ mapping to a barycenter of $\Delta_B$,
then
there exists a modified barycentric subdivision $MBS(\Delta_X)$ having
multiplicity strictly less than the multiplicity of $\Delta_X$, such
that $f$
induces a simplicial map $f':MBS(\Delta_X)\rar BS(\Delta_B)$. 

For every singular cone $\sigma\in\Delta_X$ choose a point $w_\sigma$ as
follows. By assumption, there exists a Waterman point $w\in\sigma$
mapping to
a barycenter of $\Delta_B$: $f(w)=\hat{\tau}$. Write
$f(\sigma)=\tau*\tau_0$
and choose a face $\sigma_0\leq\sigma$ such that
$f:\sigma_0\stackrel{\simeq}{\rar}\tau_0$. Set
$w_\sigma=w+\hat{\sigma}_0$;
then 
\[ f(w_\sigma) = f(w)+f(\hat{\sigma}_0)=\hat{\tau}+\hat{\tau}_0 =
\widehat{f}(\sigma)\]
Having chosen the set 
$\{w_\sigma\}$, we may remove some of the points $w_\sigma$ if necessary
so
that every simplex $\rho\in\Delta_X$ contains at most one $w_\sigma$ in
its
interior. With $\bar{\Delta}_X$ as in Example~\ref{example}, let
$\tilde{\Delta}_X =\bar{\Delta}_X \cup \{\rho\in\Delta_X|
w_\sigma\in\mbox{int $(\rho)$ for some singular $\sigma$}\}$, $b_\rho =
\hat{\rho}$ if $\rho\in\bar{\Delta}_X$, and $b_\rho = w_\sigma$ if
$w_\sigma\in\mbox{int}(\rho)$. 

By construction, $(\tilde{\Delta}_X,\{b_\rho\})$ satisfies the
hypothesis of
Corollary~\ref{cor-simpl}, hence $f$ induces a simplicial map
$f':MBS(\Delta_X)\rar BS(\Delta_B)$. Before we compute the multiplicity
of
$MBS(\Delta_X)$, we choose a particular total order $\prec$ on
$\Delta_X$. Extend $\leq$ on $\Delta_X$ to a partial order $\prec_0$ by
declaring that $\sigma_1\prec_0\sigma_2$ for all (nonsingular)
$\sigma_1\in\bar{\Delta}_X$ and singular $\sigma_2\in\Delta_X$. Let
$\prec$
be an extension of $\prec_0$ 
to a total order on $\Delta_X$. With such $\prec$, if
$\sigma\in\Delta_X$ is
singular, then $b_{\tilde{\sigma}}$ is one of the points $w_\rho$.

As in Example~\ref{example}, the multiplicity of $MBS(\Delta_X)$ is not
bigger than the multiplicity of $\Delta_X$. If $\sigma\in\Delta_X$ is
singular we show by induction on the dimension of $\sigma$ that the
multiplicity of $MBS(\sigma)$ is strictly less than the 
multiplicity of $\sigma$. Let $v_1,\ldots,v_N$ be the primitive points
of
$\sigma$, and consider the
cone $\sigma'=\LA b_{\tilde{\sigma}},v_1,\ldots,v_{N-1}\RA$ in the star
subdivision of $\sigma$ at $b_{\tilde{\sigma}}=\sum_i a_i v_i$. To show
that
every maximal cone of $MBS(\sigma)$ contained in $\sigma'$ has
multiplicity
less than the multiplicity of $\sigma$, we have three cases: 
\begin{enumerate}
\item If $a_N$ = 0, then $\sigma'$ is degenerate.
\item If $0<a_N<1$, then the multiplicity of 
${\LA}b_{\tilde{\sigma}},v_1,\ldots,v_{N-1}{\RA}$ is less than the
multiplicity of
$\sigma$, and since all $b_\rho=\sum_i c_i v_i$ have coefficients $0\leq
c_i
\leq 1$, further subdivisions at $b_\rho$ do not increase the
multiplicity of
${\LA}b_{\tilde{\sigma}},v_1,\ldots,v_{N-1}{\RA}$. 
\item If $a_N=1$, then $b_{\tilde{\sigma}}=w+\hat{\rho}$ for some
$\rho\leq\sigma$ and 
$w\in{\LA}v_1,\ldots,v_{N-1}{\RA}$ a Waterman point. Hence
${\LA}v_1,\ldots,v_{N-1}{\RA}$ is singular and, by induction, every
maximal
cone in $MBS({\LA}v_1,\ldots,v_{N-1}{\RA})$ has 
multiplicity less than the multiplicity of
${\LA}v_1,\ldots,v_{N-1}{\RA}$. Then also every maximal cone in
$\bfr_+ b_{\tilde{\sigma}}*MBS({\LA}v_1,\ldots,v_{N-1}{\RA})$ has
multiplicity
less than the multiplicity of
${\LA}b_{\tilde{\sigma}},v_1,\ldots,v_{N-1}{\RA}$.
\end{enumerate}

\section{Families of surfaces and 3-folds.}

{\bf Proof of Theorem~\ref{main-thm}.} It is not difficult to subdivide
$\Delta_X$ and $\Delta_B$ so that $\Delta_X$ is simplicial,
$\Delta_B$ is nonsingular, and  $f:\Delta_X\rar\Delta_B$ is a simplicial
map
(e.g. Proposition~4.4 and the remark following it in \cite{ak}).
Applying the
join 
construction we can make $f|_{\Delta_{X,i}}$ semistable without
increasing
the multiplicity of $\Delta_X$. We will show below that 
every singular simplex of $\Delta_X$ contains a Waterman point mapping
to a
barycenter of $\Delta_B$. By the previous section, there exist a
modified
barycentric subdivision and a simplicial map $f':MBS(\Delta_X)\rar
BS(\Delta_B)$, with multiplicity of $MBS(\Delta_X)$ strictly less than
the
multiplicity of $\Delta_X$. Since $f'$ is simplicial and $BS(\Delta_B)$
nonsingular, the proof is completed by induction. 

Restrict $f$ to a singular simplex $f:\sigma\rar\tau$, where $\sigma$
has
primitive points $v_{i j}, i=1,\ldots,n, j=1,\ldots,J_i$, $\tau$ has
primitive points $u_1\ldots,u_n$, and $f(v_{i j})=u_i$. Since $\sigma$
is
singular, it contains a Waterman point 
\[ w=\sum_{i,j} \alpha_{i j} v_{i j}, \qquad 0\leq\alpha_{i j}<1, \]
where not all $\alpha_{i j}=0$. Restricting to a face of $\sigma$ if
necessary we may 
assume that $w$ lies in the interior of $\sigma$, hence $0<\alpha_{i
j}$. Since $f(w)\in N_\tau$, it
follows that $\sum_j \alpha_{ij}\in\bfz$ for all $i$. In particular, if
$J_{i_0}=1$ for some $i_0$ then $\alpha_{i_0 1}=0$, and $w$ lies in a
face of
$\sigma$. So  we may assume that $J_i>1$ for all $i$. Since the relative
dimension of $f$ is $\sum_i (J_i-1)$, we have to consider all possible
decompositions $\sum_i (J_i-1) \leq 3$, where $J_i>1$ for all $i$.

The cases when the relative dimension of $f$ is 0 or 1 are trivial and
left
to the reader.

If the relative dimension of $f$ is 2, then either $J_1=3$,
or $J_1=J_2=2$. In the first case, we have that
${\LA}v_{11},v_{12},v_{13}{\RA}$ is singular, contradicting the
semistability of
$f|_{\Delta_{X,1}}$. In the second case, $\alpha_{11}+\alpha_{12},
\alpha_{21}+\alpha_{22} \in \bfz$ and $0< \alpha_{i j} < 1$ imply that
$\alpha_{11}+\alpha_{12}= \alpha_{21}+\alpha_{22}=1$. Hence
$f(w)=u_1+u_2$ is
a barycenter. 

In relative dimension 3, either $J_1=4$, or $J_1=3,J_2=2$, or
$J_1=J_2=J_3=2$. In
the first case, we get a contradiction with the semistability of
$f|_{\Delta_{X,1}}$; the third case gives $\alpha_{11}+\alpha_{12}=
\alpha_{21}+\alpha_{22}=\alpha_{31}+\alpha_{32}=1$ as for relative
dimension
$2$. In the second case either $\alpha_{11}+\alpha_{12}+\alpha_{13} =
\alpha_{21}+\alpha_{22}=1$ and $w$ maps to a barycenter, or
$\alpha_{11}+\alpha_{12}+\alpha_{13} = 2, 
\alpha_{21}+\alpha_{22}=1$ and $(\sum v_{i j})-w$ maps to a barycenter.
\QED

\end{document}